\documentclass[11pt,letterpaper]{amsart}
\usepackage[utf8]{inputenc}
\usepackage{amsmath}
\usepackage{amsfonts}
\usepackage{amssymb}
\usepackage{mathrsfs}
\usepackage{graphicx}
\usepackage[left=2.9cm,right=2.9cm, top=2.9cm,bottom=2.9cm]{geometry}
\usepackage[colorlinks=true, linkcolor=black, citecolor=black, urlcolor=black]{hyperref}
 \numberwithin{equation}{subsection}

\usepackage{tikz-cd} 
\usepackage[all,cmtip]{xy}	
\xyoption{arrow}
\usepackage{xcolor,etoolbox}

\usepackage{enumitem}

\theoremstyle{definition}
\newtheorem{defi}[subsection]{Definition}

\newtheorem{rmk}[subsection]{Remark}
\newtheorem{eg}[subsection]{Example}

\theoremstyle{plain}
\newtheorem{theorem}[subsection]{Theorem}
\newtheorem{lem}[subsection]{Lemma}
\newtheorem{prop}[subsection]{Proposition}

\newtheorem*{theoremint}{Theorem}

\newcommand{\ZZ}{\mathbb{Z}}

\newcommand{\Pp}{\mathcal{P}}
\newcommand{\Ff}{{\mathcal{F}}}
\newcommand{\I}{\mathscr{I}\!}

\newcommand{\Rr}{\mathcal{R}}

\newcommand{\Xt}{{\widetilde{X}}}

\newcommand{\Hom}{\textrm{Hom}}
\newcommand{\Isom}{\mathscr{I}\!\textit{so}}
\newcommand{\Aut}{\text{Aut}}

\newcommand{\Der}{\textrm{Der}}
\newcommand{\Spec}{\textrm{Spec}}

\newcommand{\Bun}{\mathcal{B}\textit{un}}

\renewcommand{\H}{\mathcal{H}}
\newcommand{\Go}{\mathcal{G}_{\mathcal{P}}}
\newcommand{\IoP}{\mathscr{I}_{\mathcal{P}}(\mathcal{P}')}

\author[Chiara Damiolini]{Chiara Damiolini}
\address{\textrm{Chiara Damiolini} \newline \indent Department of Mathematics, University of Pennsylvania, Philadelphia, PA 19104}
\email{chiara.damiolinid@gmail.com}

\title{On equivariant bundles and their moduli spaces}
\begin{document}

\begin{abstract} 

Let $G$ be an algebraic group and $\Gamma$ a finite subgroup of automorphisms of $G$. Fix also a possibly ramified $\Gamma$-covering $\Xt \to X$. In this setting one may define the notion of $(\Gamma,G)$-bundles over $\Xt$ and, in this paper, we give a description of these objects in terms of $\H$-bundles on $X$, for an appropriate group $\H$ over $X$ which depends on the local type of the $(\Gamma,G)$-bundles we intend to parametrize. This extends, and along the way clarifies, an earlier work of Balaji and Seshadri. 
\end{abstract}

\maketitle

\section{Introduction}

The moduli space $\Bun_{\text{SL}_r}(X)$ parametrizing vector bundles of rank $r$ and trivial determinant over a smooth curve $X$ is a central object in algebraic geometry, with deep connections to representation theory and conformal field theory \cite{ramanan:1973:moduli, BeauvilleLaszlo1994Conformal, beauville:1993:CFT, sorger1996formule}. One can extend this notion by replacing $\text{SL}_r$ with another simple algebraic group $G$ over the base field or, if there are marked points on the curve, by considering \textit{parabolic} $G$-bundles \cite{KumarNR1994Grassmannian, mehta.seshadri:1980:moduli}. In all these cases there are connections to conformal field theory through generalized theta functions and conformal blocks   \cite{Pauly1996Parabolic, LS1997PicardBunG}. 
An effective way to describe all the previous instances is to use \textit{parahoric} bundles, that is principal bundles on $X$ with respect to a parahoric Bruhat-Tits group $\H$ defined over the curve $X$ itself \cite{BruhatTitsII, heinloth2010uniformization}. A natural way to produce such groups, and which  gives rise to every \textit{split} parahoric Bruhat-Tits group, is from Galois coverings of curves \cite{balaji2011moduli}. In this circumstance as well, links to conformal field theory through twisted 
conformal blocks have been established \cite{damiolini:2020:conformal,zelaci:2019:moduli,hong2021conformal}.  

We explain how groups can be constructed from coverings. Let $\Gamma$ be a finite group and let $\pi \colon \Xt \to X$ be a possibly ramified $\Gamma$-covering of smooth curves. Assume that $G$ is an algebraic group which is equipped with a group homomorphism $\rho \colon \Gamma \to \Aut(G)$. Then the group $\pi_*(\Xt \times G)$ has a natural action of $\Gamma$, and so one may construct the subgroup of $\Gamma$-invariant elements $\H:= (\pi_*(\Xt \times G))^\Gamma$, which is a parahoric Bruhat-Tits group over $X$.  The construction of a group starting with a covering as above still holds without assuming that $\Xt$ and $X$ are curves. In this paper we aim to characterize $\H$-bundles in terms of certain $G$-bundles on $\Xt$. A particular instance of our main result (Theorem \ref{thm-equivstacks}) can be stated as follows:

\begin{theoremint} The functor $\pi_*(-)^\Gamma$ induces an equivalence between $\Bun^{\text{triv}}_{(\Gamma,G)}(\Xt)$ and $\Bun_\H(X)$. \end{theoremint}

Here $\Bun_{(\Gamma,G)}(\Xt)$ is the stack parametrizing $(\Gamma, G)$-bundles on $\Xt$, that is $G$-bundles on $\Xt$ which are equipped with an induced action of $\Gamma$ which lifts the action of $\Gamma$ on $\Xt$ and is compatible with the action of $\Gamma$ on $G$ (see Definition \ref{defi:GammaGbundle}). Given every $(\Gamma,G)$-bundle $\Pp$, the $\Bun^{\Pp}_{(\Gamma,G)}(\Xt)$ is given by those $(\Gamma,G)$-bundles which have the same \textit{local type} as $\Pp$. When $\Pp$ is the trivial $(\Gamma, G)$-bundle, then we use the notation  $\Bun^{\text{triv}}_{(\Gamma,G)}(\Xt)$ instead. A fundamental ingredient to understand our main result is the concept of \textit{local types} of $(\Gamma,G)$-bundles (Definition \ref{defi:localtype}). Essentially two $(\Gamma,G)$-bundles have the same local type exactly when they are locally isomorphic not only as $G$-bundles, but as $(\Gamma, G)$-bundles, i.e. taking into account also the action of $\Gamma$. 

From another point of view, $\Bun^{\text{triv}}_{(\Gamma,G)}(\Xt)$ can be also interpreted as the biggest substack of $\Bun_{(\Gamma,G)}(\Xt)$ where we can apply the functor $\pi_*(-)^\Gamma$ and land in $\Bun_\H$. In Example \ref{eg:S4} we show that $\Bun^{\text{triv}}_{(\Gamma,G)}(\Xt)$ is a strict substack of $\Bun_{(\Gamma,G)}(\Xt)$, emphasizing in this way the fact that not all $(\Gamma,G)$-bundles are locally trivial, and that fixing local types is a necessary condition for our theorem to hold. 

\subsection*{Comparison with \cite{balaji2011moduli}} The inspiration for working with bundles associated with groups arising from coverings, and to give their description in terms of $(\Gamma,G)$-bundles comes from Balaji's and Seshadri's paper \cite{balaji2011moduli}, where the authors give a description of the moduli space of parahoric torsors over a smooth curve. In particular, in \cite[Theorem 4.1.6]{balaji2011moduli} a similar statement to our main result is asserted with two main differences.

The first, and perhaps most important, concerns the concept of local type. More precisely Theorem 4.1.6 states that $\pi_*(-)^\Gamma$ is an equivalence between $\Bun_{(\Gamma,G)}(\Xt)$ and $\Bun_\H(X)$. However, as illustrated by Example \ref{eg:S4}, this is not true in general. As shown in Proposition \ref{prop-PisHbundle}, in order to address this problem, we need to introduce and use the concept of local types. 

The second difference lies in how $\Gamma$ acts on $G$: in \cite{balaji2011moduli} the authors let $\Gamma$ act on $G$ via inner automorphisms only. Under this hypothesis, and assuming that $\Xt$ and $X$ are smooth curves, $\H = \pi_*(-)^\Gamma$ is a \textit{split} parahoric Bruhat-Tits group. We do not impose any condition on how $\Gamma$ acts on $G$, obtaining in this case a wider class of groups $\H$ arising from coverings.

In conclusion, from the above considerations, it follows that our main result not only extends \cite[Theorem 4.1.6]{balaji2011moduli} beyond the case of inner automorphisms, but can be seen also as a way to clarify the need of fixing local types of $(\Gamma,G)$-bundles for that statement to hold true.

\subsection*{Further projects} As aforementioned, one of the main motivations of this work comes from the study of parahoric Bruhat-Tits groups over a curve, so we will assume in this section that $\Xt$ and $X$ are smooth curves. In this setting, \cite[Theorem 5.2.7]{balaji2011moduli} states that all \textit{split} parahoric Bruhat-Tits can be recovered from coverings, provided that the Galois group $\Gamma$ acts on $G$ via inner automorphisms only. However 
 it is natural to include also automorphism which are not necessarily inner (see for instance \cite[Example 2.3]{damiolini:2020:conformal}). It follows that the category of groups $\H$ which can be constructed through coverings, and to which we can apply our main result, is much larger than the one studied in \cite{balaji2011moduli}. We can then ask whether all parahoric Bruhat-Tits groups arise in this fashion. If the answer to this question is affirmative, then Theorem \ref{thm-equivstacks} implies that the study of all parahoric bundles over a curve $X$ can be translated into the study of $\Gamma$-equivariant bundles over $\Xt$ for an appropriate $\Gamma$-covering $\Xt \to X$. Observe that if this were true, then \cite{damiolini:2020:conformal, hong2021conformal} would provide the notion of conformal blocks associated to every parahoric Bruhat-Tits group.

Another element under investigation is the classification of the possible local types. In the case of inner automorphisms, one can describe local types via conjugation classes in $G$. However an explicit classification of local types for more general automorphisms has not been given, and it is one of the central topic of a forthcoming paper of the author together with Jiuzu Hong \cite{damiolini.hong}. 

\subsection*{Plan of the paper} We begin the paper by fixing some notation which will be used throughout. We then introduce the main ingredients of the paper: we explain what we mean by $(\Gamma,G)$-bundles in Definition \ref{defi:GammaGbundle} and, after Example \ref{eg:S4}, we introduce the notion of local type in Definition \ref{defi:localtype}. This leads to Proposition \ref{prop-PisHbundle} which identifies the correct substack of $\Bun_{(\Gamma,G)}$ to which we can apply the functor $\pi_*(-)^\Gamma$ and obtain $\H$-bundles. Finally, in Section \ref{sec:mainthm} we state and prove our main result, Theorem \ref{thm-equivstacks} by explicitly constructing the map providing the inverse to $\pi_*(-)^\Gamma$ (Proposition \ref{prop-inversedef}).

\section{$(\Gamma,G)$-bundles and local types} \label{sec:loctype}

Throughout this paper, we will fix a finite group $\Gamma$. All the schemes will be over a field $k$ whose characteristic does not divide the order of $\Gamma$. We will further fix the following data:
\begin{itemize}
\item An affine algebraic group $G$ over $k$ which is endowed with a group homomorphism $\rho \colon \Gamma \to \Aut(G)$ sending $\gamma$ to $\gamma_G$.
\item A (ramified) $\Gamma$-covering $\pi \colon \Xt \to X$  over $k$, that is
\begin{itemize} \item $\pi$ is a finite flat morphism; \item the group of automorphisms of $\Xt$ over $X$ is isomorphic to $\Gamma$. The automorphism of $\Xt$ associated with $\gamma \in \Gamma$ is denoted $\gamma_{\Xt}$; \item $\Xt$ is a generically étale $\Gamma$-torsor over $X$ via $\pi$.\end{itemize}
\end{itemize}



The \emph{ramification locus} of $\pi$ is the subscheme of $\Xt$ where $\pi$ is not étale.  Its image in $X$ is denoted by $\Rr$ and called the \emph{branch locus} of $\pi$.

\begin{defi} Given a $G$-bundle $\Pp$ on $\Xt$ we denote by $\Go$ the automorphisms group scheme $\Isom_G(\Pp,\Pp)$. 
For any $G$-bundle $\Pp'$ on $\Xt$ the scheme $\IoP:=\Isom_G(\Pp,\Pp')$ is a $\Go$-bundle.\end{defi}
 Recall (see for instance \cite[Section 7.6]{bosch1990neron}) that the Weil restriction $\pi_*\IoP$ of $\IoP$ along $\pi$ is defined by the equality \[\pi_*\IoP(T) := \Hom_\Xt(T \times_X \Xt, \IoP)\] for every $T$ over $X$. It follows from \cite[Section 7.6]{bosch1990neron} that $\pi_*\IoP$ is representable by a smooth scheme over $X$ and that $\pi_*\Go$ has the structure of an algebraic group. 
\medskip

\noindent The following statement is a version of \cite[Lemma 4.1.4]{balaji2011moduli}, which we add for completeness. 

\begin{lem} \label{lem-bundle} Let $\Pp'$ be a $G$-bundle over $\Xt$, then $\pi_*\IoP$ is a $\pi_*\Go$-bundle.
\end{lem}

\proof Since taking fibred products and Weil restrictions commute, $\pi_*\Go$ still acts on $\pi_*\IoP$. Similarly we have that $\pi_*\IoP \times_X \pi_*\Go \cong \pi_*\IoP \times_X \pi_*\IoP$ via the map $(f,g) \mapsto (f, fg)$, thus we are left to prove that for every point $x \in X(\overline{k})$ there exists an étale neighbourhood $U$ of $x$ such that $(\pi_*\IoP)(U) \neq \emptyset$. Since $\pi$ is finite we know that $\pi^{-1}\lbrace x \rbrace$ is a finite scheme over $\Spec(\overline{k})$ over which both $\Pp$ and $\Pp'$ are trivial. It follows that the map $q \colon \pi_*\IoP \to X$ is surjective. We conclude that $q$ is smooth and surjective. Applying  \cite[Corollaire 17.16.3]{EGAIV4} for every $x \in X$ there exists an étale neighbourhood $U$ of $x$ such that $(\pi_*\IoP)(U) \neq \emptyset$.\endproof

We will be interested in those $G$-bundles over $\Xt$ which are equipped with an action of $\Gamma$, compatible with the action of $\Gamma$ on $\Xt$.

\begin{defi} \label{defi:GammaGbundle} A \emph{$(\Gamma, G)$-bundle} on $\Xt$ is a right $G$-bundle $\Pp$ together with a left action of $\Gamma$ on its total space lifting the action of $\Gamma$ on $\Xt$ and which is compatible with the action of $\Gamma$ on $G$. The automorphism of $\Pp$ lifting $\gamma_\Xt$ will be denoted $\gamma_\Pp$. The compatibility with the action of $\Gamma$ on $G$ means that the equality
 \begin{equation} \label{eq:GammaG} \gamma_{\Pp}(p g)=\gamma_{\Pp}(p) \, \gamma_G(g)
\end{equation} holds for all $p \in \Pp$ and $g \in G$. We say that a $(\Gamma,G)$-bundle $\Pp$ on $\Xt$ is \textit{trivial} if it is isomorphic to the trivial $G$-bundle and $\gamma_\Pp = \gamma_G$ for every $\gamma \in \Gamma$.
\end{defi}

The goal of this paper is to describe $(\Gamma,G)$-bundles over $\Xt$ in terms of $\H$ bundles over $X$ for an appropriate group scheme $\H$ over $X$. We show how to attach to every $(\Gamma, G)$ bundle $\Pp$, a group scheme $\H_\Pp$ on $X$. Define the action of $\Gamma$ on $\Go$ lifting the action of $\Gamma$ on $\Xt$ via the map $\Gamma \to \Aut(\Go)$ sending the element $\gamma$ to the automorphism $\gamma_{\Go}$ defined by
\begin{equation*} \gamma_{\Go} (\phi) := \gamma_\Pp \circ \phi \circ  \gamma_\Pp^{-1}
\end{equation*}
for all $\gamma \in \Gamma$ and $\phi \in \Go$. Moreover, the left actions of $\Gamma$ on $\Go$ and on $\Xt$ induce a left action of $\Gamma$ on the Weil restriction $\pi_*\Go$ given by
\[ (\gamma f)(t, x) :=\gamma_{\Go} \, f(t,\gamma_\Xt^{-1}(x)) 
\] for every $(t,\widetilde{x}) \in T \times_X \Xt$ and $f \in \Hom(T \times_X \Xt, \Go)$.  We then define $\H_\Pp$ as the subgroup of $\Gamma$-invariant elements of $\pi_*(\Go)$, that is $\H_\Pp:=(\pi_*(\Go))^\Gamma$. By \cite[Proposition 3.4]{Edixhoven1992Neron}, $\H_\Pp$ is a smooth group over $X$. 

We now want to identify which $(\Gamma,G)$-bundles on $\Xt$ can be described in terms of $\H_\Pp$ bundles over $X$ and, conversely, whether $\H_\Pp$-bundles on $X$ determine a $(\Gamma,G)$-bundle on $\Xt$. 

For any $(\Gamma, G)$-bundle $\Pp'$, the scheme $\IoP$ is a $(\Gamma, \Go)$-bundle where the action of $\Gamma$ is given by
\[ (\gamma, \phi) \mapsto \gamma_{\Pp'} \circ  \phi \circ \gamma_{\Pp}^{-1}\]
for all $\gamma \in \Gamma$ and $\phi \in \IoP$. 
Moreover, as for $\pi_*\Go$, also $\pi_*\IoP$ is equipped with a left action of $\Gamma$, hence it makes sense to consider the subsheaf of $\Gamma$-invariant elements $(\pi_*\IoP)^\Gamma$.

As suggested in \cite[Lemma 4.1.5]{balaji2011moduli}, the candidate $\H_\Pp$-bundle corresponding to the $(\Gamma,G)$-bundle $\Pp'$ should be $(\pi_*\IoP)^\Gamma$. However, as we see in Example \ref{eg:S4}, this is not always the case. In fact, the group schemes $\H_\Pp$ heavily depend on how $\Gamma$ acts on $\Pp$ and not merely on the action of $\Gamma$ on $G$ given by $\rho$.

\begin{eg}\label{eg:S4} Let $\Gamma=\ZZ/2\ZZ=\{ 1, \gamma\}$ and $G=\mathfrak{S}_4$, the symmetric group on four elements which acts on $G$ via $\gamma_G(\alpha)=(34)(12)\alpha(12)(34)$. Consider the $\Gamma$-covering given by $\Xt=\Spec(k[t]) \to \Spec(k[t^2])=X$ and let $x \in X$ be the only ramification point of the covering. Let $\widetilde{G}=G \times \Xt$ be the trivial $(\Gamma,G)$-bundle on $\Xt$. Let $\Pp'$ be the $(\Gamma, G)$-bundle on $\Xt$ which is trivial as a $G$-bundle, but with the action of $\Gamma$ given by $\gamma_{\Pp'}(\alpha)=(34)(\alpha) (12)(34)$ for every $\alpha \in \mathfrak{S}_4$. Observe that both $\Pp$ and $\Pp'$ are $(\Gamma,G)$-bundles on $\Xt$ since equation \eqref{eq:GammaG} is satisfied. It follows that  
\[ (\pi_*\widetilde{G})^\Gamma(x)=\{\alpha \in \mathfrak{S}_4 \, | \, \alpha=(34)(12)\alpha(12)(34)\} \neq \emptyset,\]
since for instance the identity is an element of that set. We can see that 
\[ (\pi_*\Pp')^\Gamma(x)=\{\alpha \in \mathfrak{S}_4 \, | \, \alpha=(34)\alpha(12)(34)\} = \emptyset\]
for parity reasons (see also Lemma \ref{lem:reduced}).
\end{eg}

The example above shows that Equation \eqref{eq:GammaG} does not automatically imply that all $(\Gamma, G)$-bundles are locally isomorphic to the trivial $(\Gamma,G)$-bundle.  

\begin{rmk} Observe that from  \cite[Lemma 4.1.5]{balaji2011moduli} it should follow that $\pi_*(\Pp')^\Gamma$ is a $\pi_*(\widetilde{G})^\Gamma$-bundle, but Example \ref{eg:S4} provides a counterexample to that statement. To correct this problem we introduce the concept of local type, generalizing the idea introduced in \cite{balaji2011moduli}.
\end{rmk}  

\begin{defi} \label{defi:localtype} Let $\Pp_1$ and $\Pp_2$ be two $(\Gamma,G)$-bundles on $\Xt$. Then they have the \emph{same local type} at $x \in X(\overline{k})$ if \[\Isom_G\left(\Pp_1 \times \pi^{-1}\{x\} , \Pp_2 \times \pi^{-1} \{x\}  \right)^\Gamma \neq \emptyset.\]
We say that $\Pp_1$ and $\Pp_2$ have the \emph{same local type}, and we write $\Pp_1 \sim \Pp_2$, if they have the same local type at any geometric point of $X$.\end{defi}

Although the two following lemmas will not be used to show the main result, we report them as a tool to the reader who wishes to compute local types of $(\Gamma,G)$ bundles in explicit cases. A more exhaustive treatment of local types will appear in \cite{damiolini.hong}.

\begin{lem} \label{lem:reduced} Two $(\Gamma, G)$-bundles $\Pp_1$ and $\Pp_2$ have the same local type if and only if \[\Isom_G\left(\Pp_1 \times (\pi^{-1}\{x\})_{red}, \Pp_2 \times (\pi^{-1}\{x\})_{red}\right)^\Gamma \neq \emptyset \] for all $x \in X(\overline{k})$ .
\end{lem}

\proof It is enough to prove that if $\Isom_G\left(\Pp_1 \times (\pi^{-1}\{x\})_{red}, \Pp_2 \times (\pi^{-1}\{x\})_{red}\right)^\Gamma$ is not empty, then $\Pp_1$ and $\Pp_2$ have the same local type at $x$. Let $\pi^{-1}\{ x \}=\Spec(A)$ where $A$ is a finite Artin $k$-algebra. Let $\mathfrak{m}$ be its maximal nilpotent ideal, so that $\left(\pi^{-1}\{ x \}\right)_{red}= \Spec(A/\mathfrak{m})$. 
Let $\Spec(B)= \Isom_G(\Pp_1, \Pp_2)$ and by assumption there exists $\varphi_{0} \colon B \to A/\mathfrak{m}$ which is $\Gamma$-invariant and makes the diagram commute:
\[\xymatrix{ && B \ar[dll]_{\varphi_{0}}\\  A/\mathfrak{m}&& \ar[ll] \ar[u]A \rlap{~.}}\]
The aim is to lift $\varphi_{0}$ to a $\Gamma$-equivariant morphism $B \to A$ and we can reduce to the case $\mathfrak{m}^2=0$. Since $B$ is smooth over $A$ we know that $\varphi_{0}$ admits a lift $\varphi \colon B \to A$. For any $\gamma \in \Gamma$ the element $\gamma(\varphi)$ is another lift of $\varphi_{0}$, so the association $\gamma \mapsto \varphi-\gamma(\varphi)$ defines a map $h \colon \Gamma \to \Der_A(B,\mathfrak{m})$. Since $h$ satisfies the cocyle condition, we have that $h$ can be seen as an element of $\text{H}^1(\Gamma, \Der_A(B, \mathfrak{m}))$, which however is zero because the characteristic of $k$ does not divide the order of $\Gamma$. This means that there exists a derivation $\partial \in \Der_A(B, \mathfrak{m})$ such that $h(\gamma)=\gamma(\partial)-\partial$ for every $\gamma \in \Gamma$. This implies that the lift $\varphi_1 :=\varphi + \partial$ is a $\Gamma$-invariant lift of $\varphi_{0}$ and concludes the proof. \endproof

\begin{lem} \label{lem:loctypeR} Let $\Pp_1$ and $\Pp_2$ be two $(\Gamma,G)$-bundles on $\Xt$. Then $\Pp_1$ and $\Pp_2$ have the same local type if and only if they have the same local type at any geometric point of the branch locus $\Rr \subset X$.
\end{lem}

\proof It is sufficient to show that any two $(\Gamma, G)$ bundles have the same local type on $X \setminus \Rr$.  As $\pi$ is étale on $U:= X \setminus \Rr$ we can chose an étale covering $V\to U$ such that $\pi^{-1}(V)=\Gamma \times V$, i.e. the covering becomes trivial over $V$. We are then left to demonstrate that any $(\Gamma,G)$-structure defined on $\Gamma \times V \times  G$ via automorphisms $\gamma_\Pp$ satisfying Equation \eqref{eq:GammaG} is isomorphic to the trivial one induced by $\gamma_G$.  
Any isomorphism between these $(\Gamma,G)$-bundles is uniquely determined by tuples $(\alpha_\gamma)_{\gamma \in \Gamma}$ of isomorphisms of $V \times G$ satisfying \[\alpha_{\gamma \sigma} \circ \gamma_G = \gamma_\Pp \circ \alpha_\sigma\] for every $\gamma,\sigma \in \Gamma$. Since 
\begin{itemize} \item $\alpha_\gamma$ is uniquely determined by the element $\alpha_\gamma(1) \in G(V)$, and 
\item every automorphism $\gamma_\Pp$ on $\Gamma \times V \times G$ satisfying \eqref{eq:GammaG} is uniquely determined by the element $\gamma_\Pp(1) \in G(V)$,\end{itemize}
we deduce that defining $\alpha_\gamma(1):=\gamma_\Pp(1)$ gives rise to the wanted isomorphism.\endproof

The following Proposition tells us that $\H_\Pp$ can only detect those $(\Gamma,G)$-bundles having the same local type as $\Pp$.

\begin{prop} \label{prop-PisHbundle} Let $\Pp$ be a $(\Gamma,G)$-bundle over $\Xt$. Then the sheaf $(\pi_*\IoP)^\Gamma$ is an $\H_{\Pp}$-bundle if and only if $\Pp'$ has the same local type as $\Pp$.
\end{prop}

\proof We have already proved in Lemma \ref{lem-bundle} that $\pi_*\IoP$ is a $\pi_*\Go$-bundle, so that we have the isomorphism
\[\Psi \; \colon \; \pi_*\IoP \times_{X} \pi_*\Go \; \cong \; \pi_*\IoP  \times_{X} \pi_*\IoP \]
induced from $\IoP \times_{\Xt} \Go \cong \IoP \times_{\Xt} \IoP$ which sends $(f,g)$ to $(f, fg)$. The group $\Gamma$ acts diagonally on both source and target of $\Psi$ and with respect to this action $\Psi$ is $\Gamma$-equivariant. Thus it induces an isomorphism
\[ \Psi^\Gamma \; \colon \; \left(\pi_*\IoP\right)^\Gamma \times_X \H_{\Pp} \; \cong \; \left(\pi_*\IoP\right)^\Gamma \times_X \left(\pi_*\IoP\right)^\Gamma.\]
In order to conclude we need to check that $\left(\pi_*\IoP\right)^\Gamma$ admits local sections if and only if $\Pp'$ has the same local type as $\Pp$. Suppose that for every point $x \in  X$ there exists an étale neighbourhood $(u,U) \to (x,X)$ of $x$ such that there exists $\phi \in \left(\pi_*\IoP\right)^\Gamma(U)$. This implies that the composition $\phi u$ is an element in $\left(\pi_*\IoP \right)^\Gamma(x)$  which means that $\Pp$ and $\Pp'$ have the same local type at $x$. Conversely, assume that $\Pp'$ and $\Pp$ have the same local type. By definition this means that $\left(\pi_*{\IoP}\right)^\Gamma(x) \neq \emptyset$ for every geometric point $x$. It follows that the map $q \colon \left( \pi_*\IoP \right)^\Gamma \to X$ is surjective on geometric points and since it is smooth, then $q$ is surjective. Invoking \cite[Corollaire 17.16.3]{EGAIV4} we can then conclude that for every $x \in X$, the map $q$ admits a section on an étale neighbourhood $U$ of $x$, and so $\left(\pi_*\IoP\right)^\Gamma(U) \neq \emptyset$. \endproof

\section{The equivalence between $\Bun_{\H_\Pp}$ and $\Bun_{(\Gamma,G)}^{\Pp}$} \label{sec:mainthm}

In this section we state and prove our main result, Theorem \ref{thm-equivstacks}. Let $\Bun_{(\Gamma,G)}^{\Pp}$ be the stack  parametrizing $(\Gamma, G)$-bundles on $\Xt$ which have the same local type as $\Pp$ and let $\Bun_{\H_\Pp}$ be the stack parametrizing $\H_\Pp$-bundles over $X$. Using this terminology, Proposition \ref{prop-PisHbundle} implies that
\[\pi_*\I_{\Pp}(-)^\Gamma  \colon \Bun_{(\Gamma,G)}^{\Pp} \to \Bun_{\H_{\Pp}}\] is a well defined map. 
We now construct the inverse following the proof of  \cite[Theorem 4.1.6]{balaji2011moduli}. Note that the inclusion of $\H_\Pp$ inside $\pi_*\Go$ induces, by adjunction, the map $\pi^*\H_\Pp \to \Go$. It follows that $\Pp$, which is naturally a left $\Go$-bundle, has an induced left action of $\pi^*\H_\Pp$. This enables us to associate to every $\H_\Pp$-bundle, the $G$-bundle $\pi^*(\Ff) \times^{\pi^*{\H_\Pp}} \Pp$ on $X$. We can say more:

\begin{prop} \label{prop-inversedef} The assignment $\Ff \mapsto \pi^*(\Ff) \times^{\pi^*{\H_\Pp}} \Pp$ defines the map \[ \pi^*(-) \times^{\pi^*{\H_\Pp}} \Pp  \colon \Bun_{\H_{\Pp}} \to \Bun_{(\Gamma,G)}^{\Pp}.\]
 \end{prop} 
 
\proof We first show that for any $\H_\Pp$-bundle $\Ff$ over $X$, the scheme $\Ff^\Pp:=\pi^*(\Ff) \times^{\pi^*{\H_\Pp}} \Pp$ is indeed a $(\Gamma, G)$-bundle. Observe that it has a natural right action of $G$ and a left action of $\Gamma$ induced by the ones on $\Pp$. Let $\gamma \in \Gamma$ and $g \in G$ and consider $(f,p) \in \Ff^{\Pp}$. The chain of equalities
\[\gamma_{\Ff^\Pp}((f,p)g))=\gamma_{\Ff^\Pp}(f,pg)=(f, \gamma_\Pp(pg))=(f,\gamma_\Pp(p) \gamma_G(g))=\left(\gamma_{\Ff^\Pp} (f,p)\right) \gamma_G(g)\] tell us that $\Ff^{\Pp}$ is a $(\Gamma, G)$-bundle on $\Xt$.

We now check that $\Ff^{\Pp}$ has the same local type as $\Pp$. Let $x$ be a geometric point of $X$. Then there is an isomorphism
\begin{align*}
\Ff^{\Pp} \times_{\widetilde{X}} \pi^{-1}\lbrace x \rbrace &= \left(\pi^*\Ff \times_{\widetilde{X}} \pi^{-1}\lbrace x \rbrace\right) \times^{\pi^*\H_\Pp \times \pi^{-1}\lbrace x \rbrace}\left(\Pp \times_{\widetilde{X}} \pi^{-1}\lbrace x \rbrace\right) \\
&= \pi^*\left(\Ff|_x\right)\times^{\pi^*(\H_\Pp|_x )}\left(\Pp\times_{\widetilde{X}} \pi^{-1}\lbrace x \rbrace\right)  \\
&\cong \pi^*\left(\H_\Pp|_x \right)\times^{\pi^*(\H_\Pp|_x)}\left(\Pp\times_{\widetilde{X}} \pi^{-1}\lbrace x \rbrace\right) \\
&= \Pp\times_{\widetilde{X}} \pi^{-1}\lbrace x \rbrace,
\end{align*}
which is $\Gamma$-invariant because it is induced by the isomorphism between $\Ff|_x$ and $\H_\Pp|_x$ on which $\Gamma$ acts trivially.\endproof

We have gathered all the ingredients to state and prove our main result.

\begin{theorem} \label{thm-equivstacks} The maps $\pi_*\I_{\Pp}(-)^\Gamma$ and $\pi^*(-) \times^{\pi^*{\H_\Pp}} \Pp$ are each other inverses, defining an equivalence between the stacks $\Bun_{(\Gamma,G)}^{\Pp}$ and $\Bun_{\H_{\Pp}}$.
\end{theorem}

\proof Let $\Ff$ be an $\H_\Pp$-bundle. We first show that $\pi_*\I_{\Pp}(\Ff^\Pp)^\Gamma$ is naturally isomorphic to $\Ff$, where we have simplified our notation, as in the proof of Proposition \ref{prop-inversedef}, by writing $\Ff^\Pp$ in place of $\pi^*(\Ff) \times^{\pi^*{\H_\Pp}} \Pp$. The assignment $f \mapsto [\phi_f \colon p \mapsto (f,p)]$ defines a morphism from $\pi^*\Ff$ to $\I_\Pp\left(\Ff^{\Pp}\right)$. By pushing it down to $X$ and taking $\Gamma$ invariants we obtain a map 
\[\Ff \to \pi_*\I_\Pp\left(\Ff_{\Pp}\right)^\Gamma
= \pi_*\I_\Pp\left(\pi^*\Ff \times^{\pi^*{\H_\Pp}}\Pp\right)^\Gamma .\]	
Locally we can check that this map is $\H_\Pp$-equivariant, and hence it must be an isomorphism since both source and target are $\H_\Pp$-bundles.

Conversely, let $\Pp'$ be a $(\Gamma, G)$-bundle with the same local type as $\Pp$. Applying $\left(\pi_*\I_\Pp(-)\right)^\Gamma$ and then $\pi^*(-)\times^{\pi^*\H_\Pp} \Pp$ to $\Pp'$ we obtain
\[\pi^*\left(\pi_*(\IoP)^\Gamma\right) \times^{\pi^*\H_\Pp} \Pp.\]
The inclusion $\left(\pi_*(\IoP\right)^\Gamma \subseteq \pi_*(\IoP)$ induces, by adjunction, the map of $\pi^*\H_\Pp$-bundles \[\pi^*\left(\pi_*(\IoP)^\Gamma\right) \to \IoP, \quad f \mapsto \alpha_f\]
which extends to a map of $G$-bundles
\[\alpha \colon \pi^*\left(\pi_*(\IoP)^\Gamma\right) \times^{\pi^*\H_\Pp} \Pp \to \IoP \times^{\pi^*\H_\Pp} \Pp, \quad (f,p) \mapsto (\alpha_f, p).\] 
The evaluation map $\beta \colon \IoP \times^{\pi^*\H_\Pp} \Pp \to \Pp'$ allows us to obtain the morphism 
\[\beta \alpha \colon \pi^*\left(\pi_*(\IoP)^\Gamma\right) \times^{\pi^*\H_\Pp} \Pp \to\Pp'\]
which we are left to show to be equivariant with respect to the actions of $\Gamma$ and $G$. Since both $\alpha$ and $\beta$ are $G$-equivariant, their composition is as well. The $\Gamma$-invariance translates to showing that $\alpha_f\left(\gamma_{\Pp}(p)\right)$ and $\gamma_{\Pp'}(\alpha_f(p))$ coincide, which holds because $\alpha_f$ is $\Gamma$-equivariant. \endproof

\subsection{Special case} When $\Pp$ is the trivial $G$-bundle on $\Xt$ with the action of $\Gamma$ given $\rho$, then we have that $\Go \cong \Xt  \times G$. If we apply Theorem \ref{thm-equivstacks} in this context, we obtain the Theorem stated in the introduction. 

\bigskip
\noindent\textit{Acknowledgments.} I am deeply indebted to V. Balaji and C. S. Seshadri for their inspirational work and perspective on parahoric bundles. Thanks to J. Heinloth, A. Gibney and J. Hong for helpful conversations.

\bibliography{BiblioGammaGbundles}

\begin{thebibliography}{KNR}

\bibitem[Bea]{beauville:1993:CFT}
A. Beauville.
\newblock Vector bundles on {R}iemann surfaces and conformal field theory.
\newblock In {\em Algebraic and geometric methods in mathematical physics
  ({K}aciveli, 1993)}, volume~19 of {\em Math. Phys. Stud.}, pages 145--166.
  Kluwer Acad. Publ., Dordrecht, 1996.

\bibitem[BL]{BeauvilleLaszlo1994Conformal}
A. Beauville and Y. Laszlo.
\newblock Conformal blocks and generalized theta functions.
\newblock {\em Comm. Math. Phys.}, 164(2):385--419, 1994.

\bibitem[BLR]{bosch1990neron}
S. Bosch, W. L\"utkebohmert, and M. Raynaud.
\newblock {\em N\'eron models}, volume~21 of {\em Ergebnisse der Mathematik und
  ihrer Grenzgebiete (3) [Results in Mathematics and Related Areas (3)]}.
\newblock Springer-Verlag, Berlin, 1990.

\bibitem[BS]{balaji2011moduli}
V. Balaji and C.~S. Seshadri.
\newblock Moduli of parahoric {$\mathscr{G}$}-torsors on a compact {R}iemann
  surface.
\newblock {\em J. Algebraic Geom.}, 24(1):1--49, 2015.

\bibitem[BT]{BruhatTitsII}
F. Bruhat and J. Tits.
\newblock Groupes r\'eductifs sur un corps local. {II}. {S}ch\'emas en groupes.
  {E}xistence d'une donn\'ee radicielle valu\'ee.
\newblock {\em Inst. Hautes \'Etudes Sci. Publ. Math.}, (60):197--376, 1984.

\bibitem[Dam]{damiolini:2020:conformal}
C. Damiolini.
\newblock Conformal blocks attached to twisted groups.
\newblock {\em Math. Z.}, 295(3-4):1643--1681, 2020.

\bibitem[DH]{damiolini.hong}
C. Damiolini and J. Hong.
\newblock Local types of $({\Gamma},{G})$-bundles.
\newblock in prep.

\bibitem[Edi]{Edixhoven1992Neron}
B. Edixhoven.
\newblock N\'eron models and tame ramification.
\newblock {\em Compositio Math.}, 81(3):291--306, 1992.

\bibitem[Gro]{EGAIV4}
A. Grothendieck.
\newblock {\'El\'ements de g\'eom\'etrie alg\'ebrique. {IV}. \'Etude locale des
  sch\'emas et des morphismes de sch\'emas {IV}}.
\newblock {\em Inst. Hautes \'Etudes Sci. Publ. Math.}, (32):361, 1967.

\bibitem[Hei]{heinloth2010uniformization}
J. Heinloth.
\newblock Uniformization of {$\mathscr{G}$}-bundles.
\newblock {\em Math. Ann.}, 347(3):499--528, 2010.

\bibitem[HK]{hong2021conformal}
J. Hong and S. Kumar.
\newblock Conformal blocks for galois covers of algebraic curves, 2021.

\bibitem[KNR]{KumarNR1994Grassmannian}
S. Kumar, M.~S. Narasimhan, and A. Ramanathan.
\newblock Infinite {G}rassmannians and moduli spaces of {$G$}-bundles.
\newblock {\em Math. Ann.}, 300(1):41--75, 1994.

\bibitem[LS]{LS1997PicardBunG}
Y. Laszlo and C. Sorger.
\newblock The line bundles on the moduli of parabolic {$G$}-bundles over curves
  and their sections.
\newblock {\em Ann. Sci. \'Ecole Norm. Sup. (4)}, 30(4):499--525, 1997.

\bibitem[MS]{mehta.seshadri:1980:moduli}
V.~B. Mehta and C.~S. Seshadri.
\newblock Moduli of vector bundles on curves with parabolic structures.
\newblock {\em Math. Ann.}, 248(3):205--239, 1980.

\bibitem[Pau]{Pauly1996Parabolic}
C. Pauly.
\newblock Espaces de modules de fibr\'es paraboliques et blocs conformes.
\newblock {\em Duke Math. J.}, 84(1):217--235, 1996.

\bibitem[Ram]{ramanan:1973:moduli}
S. Ramanan.
\newblock The moduli spaces of vector bundles over an algebraic curve.
\newblock {\em Math. Ann.}, 200:69--84, 1973.

\bibitem[Sor]{sorger1996formule}
C. Sorger.
\newblock La formule de {V}erlinde.
\newblock {\em Ast\'erisque}, (237):Exp.\ No.\ 794, 3, 87--114, 1996.
\newblock S\'eminaire Bourbaki, Vol. 1994/95.

\bibitem[Zel]{zelaci:2019:moduli}
H. Zelaci.
\newblock Moduli spaces of anti-invariant vector bundles and twisted conformal
  blocks.
\newblock {\em Math. Res. Lett.}, 26(6):1849--1875, 2019.

\end{thebibliography}
\bibliographystyle{alphanumN}
\vfill 

\end{document}